\documentclass[12pt,reqno]{article}
\usepackage{amsmath,amssymb,amsfonts}
\usepackage{mathrsfs}
\usepackage{multicol,multirow}
\usepackage{rotating}
\usepackage{appendix}
\usepackage[numbers]{natbib}
\usepackage{tabstackengine}
\usepackage{graphics}
\usepackage{fullpage}
\usepackage{xcolor}
\usepackage{mathtools}
\usepackage{scalerel}
\usepackage{xparse}
\usepackage{amsthm}
\usepackage{amscd}
\usepackage{stmaryrd}
\usepackage{latexsym}
\usepackage{epsf}
\usepackage{breakurl}
\usepackage{hyperref}
\usepackage{color}
\usepackage{caption}
\usepackage[english]{babel}
\usepackage{diagbox}
\usepackage{stmaryrd}
\usepackage{nicematrix}
\usepackage{tabstackengine}
\usepackage{hyperref}
\RequirePackage{luatex85}
\usepackage[all]{xy}
\newtheorem{thm}{Theorem}[section]
\newtheorem{cor}[thm]{Corollary}
\newtheorem{lem}[thm]{Lemma}

\newtheorem{prop}[thm]{Proposition}
\theoremstyle{definition}
\newtheorem{defn}[thm]{Definition}
\theoremstyle{remark}
\newtheorem{rem}[thm]{Remark}
\newtheorem*{ex}{Example}
\numberwithin{equation}{section}
\DeclareMathOperator{\bl}{bl}

\newcommand{\lah}[2]{\genfrac \lfloor \rfloor {0pt} {0} {#1} {#2}} 
\newcommand{\laht}[2]{\genfrac \lfloor \rfloor {0pt} {1} {#1} {#2}}

\newcommand{\sttwo}[2]{\genfrac \{ \} {0pt} {0} {#1} {#2}}
\newcommand{\lrlah}[4]{{\genfrac \lfloor \rfloor {0pt} {0} {#1} {#2}}^{(#3)}_{#4}} 
\newcommand{\lrlaht}[4]{{\genfrac \lfloor \rfloor {0pt} {1} {#1} {#2}}^{(#3)}_{#4}} 

\begin{document}

	\title{The Lah Numbers with Higher Level and the Lah Numbers of Order $s$}
	
	\author{Aleks \v{Z}igon Tankosi\v{c} \\ 
		Gimnazija Nova Gorica \\
		Delpinova ulica 9\\
		5000 Nova Gorica \\
		Slovenia \\
		\url{zigontaleks@gmail.com}}
	
	\date{October 30, 2025}
	
	\maketitle

	\begin{abstract}
		In this paper we introduce and study two generalizations of Lah numbers, analogous to the Stirling numbers with higher level - a combinatorial one (Lah numbers with higher level) and an algebraic one (Lah numbers of order $s$). We define the Lah numbers with higher level following a combinatorial approach and the Lah numbers of order $s$ following an algebraic approach. We prove a direct connection between the Lah numbers with higher level and the $(l,r)$-Lah numbers. Some properties of the Lah numbers of order $s$ and Lah polynomials of order $s$ are given. Finally, we prove connections between these two generalizations.
	\end{abstract}
	


	\section{Introduction}

	The \textit{Lah numbers} $\laht{n}{k}$ were introduced by the Slovenian mathematician Ivo Lah \cite{lah1954} in 1954 and they count the number of partitions of the set $[n] := \{1,2,\dots,n\}$ into $k$ non-empty linearly ordered blocks (\textit{lists}, for short). The combinatorial interpretation was given in \cite{petkovvsek2007combinatorial} and \cite{petkovvsek2002lahidentity}.
	
	For $n, k \ge 1$, they satisfy the recurrence 
	\begin{align*}
		\label{LahRec} \lah{n}{k} = (n+k-1) \lah{n-1}{k} +
		\lah{n-1}{k-1}
	\end{align*} which, together with the initial conditions $\laht{0}{0} = 1$,  $\laht{n}{0} = \laht{0}{k} = 0$ for $n, k > 0$, can be used to compute $\laht{n}{k}$ for all $n, k \ge 0$. Another interesting property of Lah numbers is that they appear in the coefficients in the expansion of rising powers 
	\[ x^{\overline{n}} = \prod_{i=0}^{n-1} (x+i) = x (x+1) (x+2) \cdots (x+n-1) \] 
	in terms of falling powers 
	\[ x^{\underline{n}} = \prod_{i=0}^{n-1} (x-i) = x(x-1)(x-2) \cdots (x-n+1) \] 
	and vice versa, namely: \[ x^{\overline{n}}\ =\ \sum_{k=0} ^ {n} \lah{n}{k} x^{\underline{k}},\qquad x^{\underline{n}}\ =\	 \sum_{k=0} ^ {n} (-1)^{n+k}\lah{n}{k} x^{\overline{k} } \] for all $n \ge 0$ and all $x$, where $x^{\overline{0}} = x^{\underline{0}} = 1$. These identities can easily be proved by induction on $n$.

	Tweedie \cite{tweddie} defined the Stirling numbers of both kinds with higher level. A combinatorial approach was given in \cite{komatsu2021a, komatsu2021}. \textit{Stirling numbers of the first kind with higher level} and \textit{Stirling numbers of the second kind with higher level} count ordered $s$-tuples $(\pi_1, \pi_2, \dots, \pi_s)$ of partitions of $[n]$ into $k$ cyclically ordered blocks, respectively into $k$ unordered blocks such that  
	\[ \bl \pi_{1} = \bl \pi_{2} = \dots = \bl \pi_{s}. \] 
	Here, the set of block leaders of a partition $\pi = \{b_1, b_2, \dots, b_k\}$ with blocks $b_1, \ldots, b_k$ is defined as 
	\[ \bl \pi = \{\min b_1, \min b_2, \dots, \min b_k\}, \] i.e., as the minima of the blocks $\pi$.

	Note that in \cite{komatsu2021}, the authors used the notation min$(\pi)$ instead of $\bl \pi$. Both notations have the same meaning.

	The Stirling numbers of the first kind with higher level are denoted by
	
	\begin{align*} 
		\begin{NiceArrayWithDelims}{\llbracket}{\rrbracket}{c}
			n \\ k 
		\end{NiceArrayWithDelims}_{s},
	\end{align*} and the Stirling numbers of the second kind with higher level are denoted by
	
	\begin{align*}
		\left\{\mkern-4mu\braceVectorstack{n\\k}\mkern-4mu\right\}_{s}.
	\end{align*} These numbers satisfy triangular recurrences
	
	\begin{align}
		\label{sstirling1rec}
		\begin{NiceArrayWithDelims}{\llbracket}{\rrbracket}{c}
			n \\ k 
		\end{NiceArrayWithDelims}_{s} =\begin{NiceArrayWithDelims}{\llbracket}{\rrbracket}{c}
			n-1 \\ k-1 
		\end{NiceArrayWithDelims}_{s} + (n-1)^{s}  \begin{NiceArrayWithDelims}{\llbracket}{\rrbracket}{c}
			n-1 \\ k 
		\end{NiceArrayWithDelims}_{s}
	\end{align} and
	
	\begin{align}
		\label{sstirling2rec}
		\left\{\mkern-4mu\braceVectorstack{n\\k}\mkern-4mu\right\}_{s} = \left\{\mkern-4mu\braceVectorstack{n-1\\k-1}\mkern-4mu\right\}_{s} + k ^{s} \left\{\mkern-4mu\braceVectorstack{n-1\\k}\mkern-4mu\right\}_{s},
	\end{align} with
	\begin{align}
		\label{cases1a}
		\begin{NiceArrayWithDelims}{\llbracket}{\rrbracket}{c}
			n \\ n
		\end{NiceArrayWithDelims}_{s} = \left\{\mkern-4mu\braceVectorstack{n\\n}\mkern-4mu\right\}_{s} = 1 
	\end{align} and 
	\begin{align}
		\label{cases1b}
		\begin{NiceArrayWithDelims}{\llbracket}{\rrbracket}{c}
			n \\ 0
		\end{NiceArrayWithDelims}_{s} = \begin{NiceArrayWithDelims}{\llbracket}{\rrbracket}{c}
			0 \\ n
		\end{NiceArrayWithDelims}_{s} = \left\{\mkern-4mu\braceVectorstack{n\\0}\mkern-4mu\right\}_{s} = \left\{\mkern-4mu\braceVectorstack{0\\n}\mkern-4mu\right\}_{s} = 0 \hspace{2mm} \text{for $n>0$}. 
	\end{align} Also, 
	\begin{align}
		\label{n<k}
		\begin{NiceArrayWithDelims}{\llbracket}{\rrbracket}{c}
			n \\ k
		\end{NiceArrayWithDelims}_{s} = \left\{\mkern-4mu\braceVectorstack{n\\ k}\mkern-4mu\right\}_{s} = 0 \hspace{2mm} \text{for $n < k$}.
	\end{align}
	
	Throughout the article, we use $\mathbb{N}$ to denote the set of all positive integers without 0 and $\mathbb{N}_{0}$ to denote the set of all positive integers including 0. 
	
	For $n, s \in \mathbb{N}$, the \textit{rising factorial with higher level} is defined as

	\[
	x_{s}^{\overline{\overline{n}}} = \prod_{i=0}^{n-1} (x+i^{s}) = x(x+1^{s})(x+2^{s}) \cdots (x+(n-1)^{s}).
	\] For $n, s \in \mathbb{N}$, the \textit{falling factorial with higher level} is defined as
	
	\[
	x_{s}^{\underline{\underline{n}}} = \prod_{i=0}^{n-1} (x-i^{s}) = x(x-1^{s})(x-2^{s}) \cdots (x-(n-1)^{s}),
	\] where $x_{s}^{\overline{\overline{0}}}$=$x_{s}^{\underline{\underline{0}}}$ = 1.

	In \cite{komatsu2021}, the authors proved that
	
	\[
	{x_{s}^{\overline{\overline{n}}}} = \sum_{k=0}^{n} \begin{NiceArrayWithDelims}{\llbracket}{\rrbracket}{c}
		n \\ k 
	\end{NiceArrayWithDelims}_{s} x ^{k}
	\] and 
	
	\[
	x^{n} = \sum_{k=0}^{n} \left\{\mkern-4mu\braceVectorstack{n\\k}\mkern-4mu\right\}_{s} {x_{s}^{\underline{\underline{n}}}}.
	\]

	In 2016, Luschny \cite[A268434, A269946]{SloaneOEIS} defined the \textit{Lah numbers of order 2 and 3} with the following recurrence relations:
	
	\begin{align*}
		T_{2}(n, k)= T_{2}(n-1, k-1) + ((n-1)^{2} + k^{2}) T_{2}(n-1, k)
	\end{align*} and
	
	\begin{align*}
		T_{3}(n, k)= T_{3}(n-1, k-1) + ((n-1)^{3} + k^{3}) T_{3}(n-1, k).
	\end{align*} In the OEIS, it is also mentioned that
	
	\begin{align}
		\label{lahwith2sum}
		T_{2}(n, k) = \sum_{j=k}^{n} \begin{NiceArrayWithDelims}{\llbracket}{\rrbracket}{c}
			n \\ j 
		\end{NiceArrayWithDelims}_{2} \left\{\mkern-4mu\braceVectorstack{j\\k}\mkern-4mu\right\}_{2}
	\end{align} and 
	
	\begin{align}
		\label{lahwith3sum}
		T_{3}(n, k) = \sum_{j=k}^{n} \begin{NiceArrayWithDelims}{\llbracket}{\rrbracket}{c}
			n \\ j 
		\end{NiceArrayWithDelims}_{3} \left\{\mkern-4mu\braceVectorstack{j\\k}\mkern-4mu\right\}_{3},
	\end{align} where Luschny named the Stirling numbers with higher level of both kinds \textit{the Stirling cycle numbers of order two and three} and \textit{the Stirling set numbers of order two and three} \cite[A269944, A269945, A269947, A269948]{SloaneOEIS}. See also \cite{Luschny}. \\\\
	
	The \textit{$(s,r)$-Lah numbers} $\lrlaht{n}{k}{s}{r}$ were introduced and defined by the author in \cite{zigonthelrlah}, see also \cite[A371081, A371259, A371277, A372208]{SloaneOEIS}. They count ordered $s$-tuples $(\pi_1, \pi_2, \dots, \pi_s)$ of partitions of $[n]$ into $k$ lists such that the numbers $1, 2, \dots, r$ are in distinct lists, and $ \bl \pi_{1} = \bl \pi_{2} = \dots = \bl \pi_{s} $.
	
	\begin{rem}
		We use $s$ in place of $l$, which was originally used in \cite{zigonthelrlah}, to facilitate a clearer comparison between the properties of the Lah numbers with higher level, the $(l,r)$-Lah numbers and the Lah numbers of order $s$. 
	\end{rem} For $n \ge k \ge r \ge s$, the $(s,r)$-Lah numbers satisfy
	the recurrence relation
	\begin{align*}
		\lrlah{n}{k}{s}{r} =  (n+k-1)^s
		\lrlah{n-1}{k}{s}{r} + \lrlah{n-1}{k-1}{s}{r},
	\end{align*}
	and they can be expressed explicitly as
	\begin{align*}
		\lrlah{n}{k}{s}{r}\ &=\!\!\!\!\! \sum_{r+1 \le j_1 < j_2 < \dots < j_{n-k} \le n}
		(2j_1-2)^s (2j_2-3)^s \cdots (2j_{n-k}-(n-k+1))^s.
	\end{align*}
	Their row polynomials (for $n \ge r \ge 1$ and $l \ge 1$)
	\begin{align*}
		L_n^{s,r}(x) &= \sum_{k=0}^n \lrlah{n}{k}{s}{r} x^k
	\end{align*}
	satisfy the recurrence relation
	\begin{align*}
		L_{n+1}^{s,r}(x)\ &=\ x\, L_{n}^{s,r}(x) + \sum_{j=0}^s \binom{s}{j} n^{s-j} \sum_{i=0}^j \sttwo{j}{i} x^i \frac{d^i}{dx^i} L_{n}^{s,r}(x),
	\end{align*} where $\sttwo{j}{i}$ are Stirling numbers of the second kind, which are equivalent to the Stirling numbers of the second kind with higher level for $s=1$.
	\\\\

	Stirling numbers with higher level satisfy both combinatorial and algebraic properties analogous to non-generalized Stirling numbers. But that is not the case with the generalized Lah numbers. The Lah numbers with higher level have a combinatorial interpretation, analogous to Lah numbers, but they do not satisfy relations between falling and rising factorials with higher level and relations with Stirling numbers with higher level. These relations are satisfied by the Lah numbers of order $s$. 
	
	In the first part of this paper, we define the Lah numbers with higher level using a combinatorial approach, give a recurrence relation that they satisfy and prove a direct connection to the $(l,r)$-Lah numbers. We study a recurrence relation for polynomials connected to the Lah numbers with higher level. In the second part of this paper, we define the Lah numbers of order $s$ by generalizing identities (\ref{lahwith2sum}) and (\ref{lahwith3sum}). We give a recurrence relation and a relation that they satisfy between the falling and rising factorials with higher level. Finally, in the third part of this paper, we give some connections between the Lah numbers with higher level and the Lah numbers of order $s$. \\\\
	
	\section{The Lah Numbers with Higher Level}
	
	A natural way to define a generalization of Lah numbers following a combinatorial approach, analogous to the Stirling numbers of both kinds with higher level, is via the following combinatorial definition.
	
	\begin{defn}
		\label{sLah_def} For $s \geq 1$, the \emph{Lah numbers with higher level (level s)} $\left\lfloor\!\!\left\lfloor\!
		\begin{matrix}
			n \\ k
		\end{matrix}
		\!\right\rfloor\!\!\right\rfloor_{s}$ count ordered $s$-tuples $(\pi_1, \pi_2, \dots, \pi_s)$ of partitions of $[n]$ into $k$ lists such that \[ \bl \pi_{1} = \bl \pi_{2} = \dots = \bl \pi_{s}. \]
	\end{defn}
	
	\begin{rem}
		It is obvious that
		\begin{align*}
			\left\lfloor\!\!\left\lfloor\!
			\begin{matrix}
				n \\ n
			\end{matrix}
			\!\right\rfloor\!\!\right\rfloor_{s} = 1 \hspace{2mm} \text{for $n \geq 0$,}
		\end{align*}
		and
		\begin{align*}
			\left\lfloor\!\!\left\lfloor\!
			\begin{matrix}
				n \\ 0
			\end{matrix}
			\!\right\rfloor\!\!\right\rfloor_{s} = \left\lfloor\!\!\left\lfloor\!
			\begin{matrix}
				0 \\ n
			\end{matrix}
			\!\right\rfloor\!\!\right\rfloor_{s} =  0 \hspace{2mm} \text{for $n$ $>$ 0.}
		\end{align*}
		Note that $\left\lfloor\!\!\left\lfloor\!
		\begin{matrix}
			n \\ k
		\end{matrix}
		\!\right\rfloor\!\!\right\rfloor_{s}=0$ for $n<k$ and $\left\lfloor\!\!\left\lfloor\!
		\begin{matrix}
			0 \\ 0
		\end{matrix}
		\!\right\rfloor\!\!\right\rfloor_{s} = 1$.
	\end{rem}
	
	\begin{ex}
		Let us compute $\left\lfloor\!\!\left\lfloor\!
		\begin{matrix}
			4 \\ 3
		\end{matrix}
		\!\right\rfloor\!\!\right\rfloor_{2}$. Here $n=4$, $k=3$, $s=2$, so we need to construct all partitions of the set $[4]=\{1,2,3,4\}$ into three nonempty disjoint lists. There are 12 such partitions: 
		\begin{align*}
			\pi_1 &= \{(1),(2),(3,4)\},\\ \pi_2 &=
			\{(1),(2),(4,3)\},\\ \pi_3 &= \{(1),(2,3),(4)\},\\
			\pi_4 &= \{(1),(3,2),(4)\},\\ \pi_5 &=
			\{(1),(2,4),(3)\},\\ \pi_6 &= \{(1),(4,2),(3)\},\\
			\pi_7 &= \{(1,3),(2),(4)\},\\ \pi_8 &=
			\{(3,1),(2),(4)\},\\ \pi_9 &= \{(1,4),(2),(3)\},\\
			\pi_{10} &= \{(4,1),(2),(3)\}, \\
			\pi_{11} &=  \{(1,2),(3),(4)\}, \\
			\pi_{12} &=  \{(2,1),(3),(4)\},
		\end{align*} and the sets of their block leaders are
		
		\begin{align*}
			\bl \pi_1 &= \bl \pi_2 = \bl \pi_5 = \bl \pi_6
			= \bl \pi_9 = \bl \pi_{10} = \{1,2,3\},\\ \bl
			\pi_3 &= \bl \pi_4 =  \bl \pi_7 =  \bl \pi_8 =
			\{1,2,4\},\\ \bl \pi_{11} &= \bl \pi_{12} = \{1,3,4\}.
		\end{align*}  Now we compute the number of ordered $s$-tuples (i.e., ordered pairs) of partitions $\pi_1,\pi_2,\dots,\pi_{12}$ such that partitions in the same pair share the same set of block leaders. As there are six partitions with the set of block leaders equal to $\{1,2,3\}$, four partitions with the set of block leaders equal to $\{1,2,4\}$ and two partitions with the set of block leaders equal to $\{1,3,4\}$, we find that 
		
		\[
		\left\lfloor\!\!\left\lfloor\!
		\begin{matrix}
			4 \\ 3
		\end{matrix}
		\!\right\rfloor\!\!\right\rfloor_{2} =\ 6^2 + 4^2 + 2^{2}\ =\ 56.  \]
	\end{ex}
	
	From this combinatorial definition, analogous to combinatorial definitions of Stirling numbers of both kinds with higher level, we can observe that the Lah numbers with higher level are actually a special case of the $(s,r)$-Lah numbers.
	
	\begin{thm}
		\label{lahlrlah}
		The Lah number with higher level coincides with the $(s, 1)$-Lah number (or, equivalently $(s, 0)$-Lah number). That is for all $s \geq 1$,
		\begin{align*}
			\left\lfloor\!\!\left\lfloor\!
			\begin{matrix}
				n \\ k
			\end{matrix}
			\!\right\rfloor\!\!\right\rfloor_{s} = \lrlah{n}{k}{s}{1} = \lrlah{n}{k}{s}{0}.
		\end{align*}
	\end{thm}
	
	\begin{proof}
		The condition that the elements $1, 2, \ldots, r$ are in distinct lists is trivially satisfied when $r=0$ or $r=1$. 
	\end{proof}
	
	\begin{cor}
		From Theorem \ref{lahlrlah}, we get the following identities for the Lah numbers with higher level:
		
		• recurrence relation (for $n \ge k \ge 1$)
		\begin{equation}
			\label{slah}
			\left\lfloor\!\!\left\lfloor\!
			\begin{matrix}
				n \\ k
			\end{matrix}
			\!\right\rfloor\!\!\right\rfloor_{s} =  \left\lfloor\!\!\left\lfloor\!
			\begin{matrix}
				n-1 \\ k-1
			\end{matrix}
			\!\right\rfloor\!\!\right\rfloor_{s} + (n+k-1)^{s} \left\lfloor\!\!\left\lfloor\!
			\begin{matrix}
				n-1 \\ k
			\end{matrix}
			\!\right\rfloor\!\!\right\rfloor_{s}
		\end{equation}
		
		(equation (21) in \cite{zigonthelrlah} for $r=1$)
		
		• explicit formula 
		\begin{align*}
			\left\lfloor\!\!\left\lfloor\!
			\begin{matrix}
				n \\ k
			\end{matrix}
			\!\right\rfloor\!\!\right\rfloor_{s}\ &=\!\!\!\!\! \sum_{1 \le j_1 < j_2 < \dots < j_{n-k} \le n}
			(2j_1-2)^s (2j_2-3)^s \cdots (2j_{n-k}-(n-k+1))^s
		\end{align*}
		
		(equation (25) in \cite{zigonthelrlah} for $r=0$)
		
		• a special case for $k=n-1$ 
		\begin{align}
			\left\lfloor\!\!\left\lfloor\!
			\begin{matrix}
				n \\ n-1
			\end{matrix}
			\!\right\rfloor\!\!\right\rfloor_{s} =\ 2^s \sum_{j=1}^{n-1} j^s \label{k=n-1(1)}
		\end{align}
		
		(equation (22) in \cite{zigonthelrlah} for $r=1$)
		
		• a special case for $k=1$ 
		\begin{align*}
			\left\lfloor\!\!\left\lfloor\!
			\begin{matrix}
				n \\ 1
			\end{matrix}
			\!\right\rfloor\!\!\right\rfloor_{s} =  \left( \lah{n}{1} \right)^{s} = (n!)^{s} 
		\end{align*}
		
		(equation (23) in \cite{zigonthelrlah} for $r=1$)
		
		• the polynomials (for $n \ge 1$ and $s \ge 1$)
		\begin{align}
			\label{Lnsx}
			L_n^{s}(x) &= \sum_{k=0}^n \left\lfloor\!\!\left\lfloor\!
			\begin{matrix}
				n \\ k
			\end{matrix}
			\!\right\rfloor\!\!\right\rfloor_{s} x^k
		\end{align}
		satisfy the recurrence relation
		\begin{align*}
			L_{n+1}^{s}(x)\ &=\ x\, L_{n}^{s}(x) + \sum_{j=0}^s \binom{s}{j} n^{s-j} \sum_{i=0}^j \sttwo{j}{i} x^i \frac{d^i}{dx^i} L_{n}^{s}(x)
		\end{align*} 
		(equation (29) in \cite{zigonthelrlah} for $r=1$).
	\end{cor}

	From the recurrence relation (\ref{slah}), we get tables of values for $s=2, 3, 4$ (see tables \ref{tableone}, \ref{tableone2}, \ref{tableone3}). The sequence of Lah numbers with higher level increases quickly. For example, $\left\lfloor\!\!\left\lfloor\!
	\begin{matrix}
		6 \\ 2
	\end{matrix}
	\!\right\rfloor\!\!\right\rfloor_{25}$ is
	
	\begin{align*}
		12793287638873373780319124433790833727169139966926481069803446896427008000.
	\end{align*}

	\begin{table}[ht]
		\begin{center}
			\begin{tabular}{|c|c|c|c|c|c|c|c|}
				\hline
				\diagbox[width=2em]{$n$}{$k$} & 0& 1 &2& 3 & 4 & 5 &6  \\
				\hline
				0 & 1 & 0 & 0 & 0 &  0 & 0 & 0  \\
				\hline
				1 & 0 & 1 & 0 & 0 & 0 & 0 & 0  \\
				\hline
				2 & 0 & 4 & 1 & 0 & 0 &  0 & 0   \\
				\hline
				3 & 0 & 36 & 20 & 1 & 0 & 0 & 0  \\
				\hline
				4 & 0 & 576 & 536 & 56 & 1 & 0 & 0  \\
				\hline
				5 & 0 & 14400 & 19872 & 3280 & 120 & 1 & 0  \\
				\hline
				6 & 0 & 518400 & 988128 & 229792 & 13000 & 220 & 1  \\
				\hline
			\end{tabular}
		\end{center}
		\caption{The Lah numbers with higher level $\left\lfloor\!\!\left\lfloor\!
			\begin{matrix}
				n \\ k
			\end{matrix}
			\!\right\rfloor\!\!\right\rfloor_{2}$ for $0 \le n, k\le 6$.}
		\label{tableone}
	\end{table}
	
	\begin{table}[ht]
		\begin{center}
			\begin{tabular}{|c|c|c|c|c|c|c|c|}
				\hline
				\diagbox[width=2em]{$n$}{$k$} & 0& 1 &2& 3 & 4 & 5 &6  \\
				\hline
				0 & 1 & 0 & 0 & 0 &  0 & 0 & 0  \\
				\hline
				1 & 0 & 1 & 0 & 0 & 0 & 0 & 0  \\
				\hline
				2 & 0 & 8 & 1 & 0 & 0 &  0 & 0   \\
				\hline
				3 & 0 & 216 & 72 & 1 & 0 & 0 & 0  \\
				\hline
				4 & 0 & 13824 & 9216 & 288 & 1 & 0 & 0  \\
				\hline
				5 & 0 & 1728000 & 2004480 & 108000 & 800 & 1 & 0  \\
				\hline
				6 & 0 & 373248000 & 689264640 & 57300480 & 691200 & 1800 & 1  \\
				\hline
			\end{tabular}
		\end{center}
		\caption{The Lah numbers with higher level $\left\lfloor\!\!\left\lfloor\!
			\begin{matrix}
				n \\ k
			\end{matrix}
			\!\right\rfloor\!\!\right\rfloor_{3}$ for $0 \le n, k\le 6$.}
		\label{tableone2}
	\end{table}
	
	\begin{table}[ht]
		\begin{center}
			\begin{tabular}{|c|c|c|c|c|c|c|c|}
				\hline
				\diagbox[width=2em]{$n$}{$k$} & 0& 1 &2& 3 & 4 & 5 &6  \\
				\hline
				0 & 1 & 0 & 0 & 0 &  0 & 0 & 0  \\
				\hline
				1 & 0 & 1 & 0 & 0 & 0 & 0 & 0  \\
				\hline
				2 & 0 & 16 & 1 & 0 & 0 &  0 & 0   \\
				\hline
				3 & 0 & 1296 & 272 & 1 & 0 & 0 & 0  \\
				\hline
				4 & 0 & 331776 & 171296 & 1568 & 1 & 0 & 0  \\
				\hline
				5 & 0 & 207360000 & 222331392 & 3936064 & 5664 & 1 & 0  \\
				\hline
				6 & 0 & 268738560000 & 534025032192 & 16344449536 & 41097568 & 15664 & 1  \\
				\hline
			\end{tabular}
		\end{center}
		\caption{The Lah numbers with higher level $\left\lfloor\!\!\left\lfloor\!
			\begin{matrix}
				n \\ k
			\end{matrix}
			\!\right\rfloor\!\!\right\rfloor_{4}$ for $0 \le n, k\le 6$.}
		\label{tableone3}
	\end{table}

	\subsection{Recurrence relation for polynomials $Q_{n}^{s}(x)$} \label{qnsxsection}
	
	In this section, we study the polynomials $Q_{n}^{s}(x) = x^{n} L_{n}^{s}(x)$, because they satisfy a simpler recurrence relation than $L_{n}^{s}(x)$, involving only one sum. 
	
	\begin{thm}
		For $n \geq 1$ and $s\geq 1$, the polynomials
		\begin{align*}
			\label{rowQ}
			Q_{n}^{s}(x) = x^{n} L_{n}^{s}(x)
		\end{align*}
		satisfy the recurrence relation 
		\begin{equation}
			\label{ddQ}
			Q_{n+1}^{s} (x) = x^2 Q_{n}^{s} (x) + xA_{n}^{s}(x),
		\end{equation}
		where the polynomials $A_{n}^{s}(x)$ are for $n \geq 1$, $s \geq 1$ defined by
		\begin{align*}
			A_{n}^{s} (x) = \sum_{k=0}^{n} (n+k)^{s} \left\lfloor\!\!\left\lfloor\!
			\begin{matrix}
				n \\ k
			\end{matrix}
			\!\right\rfloor\!\!\right\rfloor_{s} x^{n+k}.
		\end{align*}
	\end{thm}
	
	\begin{proof}
		Applying formulas (\ref{slah}) and (\ref{Lnsx}), we get 
		\begin{align*}
			Q_{n+1}^{s}(x) &= x^{n+1} L_{n+1}^{s}(x) \\
			&= \sum_{k=0}^{n+1} \left\lfloor\!\!\left\lfloor\!
			\begin{matrix}
				n+1 \\ k
			\end{matrix}
			\!\right\rfloor\!\!\right\rfloor_{s} x^{n+k+1} \\
			&= \left\lfloor\!\!\left\lfloor\!
			\begin{matrix}
				n+1 \\ 0
			\end{matrix}
			\!\right\rfloor\!\!\right\rfloor_{s}  x^{n+1} + \sum_{k=1}^{n} \left\lfloor\!\!\left\lfloor\!
			\begin{matrix}
				n+1 \\ k
			\end{matrix}
			\!\right\rfloor\!\!\right\rfloor_{s} x^{n+k+1} + \left\lfloor\!\!\left\lfloor\!
			\begin{matrix}
				n+1 \\ n+1
			\end{matrix}
			\!\right\rfloor\!\!\right\rfloor_{s} x^{2n+2} \\
			&= x^{n+1} \delta_{n+1, 0} + \sum_{k=1}^{n} \left(\left\lfloor\!\!\left\lfloor\!
			\begin{matrix}
				n \\ k-1
			\end{matrix}
			\!\right\rfloor\!\!\right\rfloor_{s} + (n+k)^{s} \left\lfloor\!\!\left\lfloor\!
			\begin{matrix}
				n \\ k
			\end{matrix}
			\!\right\rfloor\!\!\right\rfloor_{s} \right) x^{n+k+1} + x^{2n+2} \\
			&=  \sum_{k=0}^{n-1} \left\lfloor\!\!\left\lfloor\!
			\begin{matrix}
				n \\ k
			\end{matrix}
			\!\right\rfloor\!\!\right\rfloor_{s} x^{n+k+2} + x^{2n+2} +  x^{n+1} \delta_{n+1, 0}  +\sum_{k=1}^{n} (n+k)^{s} \left\lfloor\!\!\left\lfloor\!
			\begin{matrix}
				n \\ k
			\end{matrix}
			\!\right\rfloor\!\!\right\rfloor_{s} x^{n+k+1} \\
			&= x^{2} \sum_{k=0}^{n}  \left\lfloor\!\!\left\lfloor\!
			\begin{matrix}
				n \\ k
			\end{matrix}
			\!\right\rfloor\!\!\right\rfloor_{s} x^{n+k} + x \sum_{k=1}^{n} (n+k)^{s} \left\lfloor\!\!\left\lfloor\!
			\begin{matrix}
				n \\ k
			\end{matrix}
			\!\right\rfloor\!\!\right\rfloor_{s} x^{n+k} \\
			&= x^{2} Q_{n}^{s}(x) + xA_{n}^{s}(x),
		\end{align*}
		where $\delta_{n,0}$ is Kronecker's delta (note that $\delta_{n+1,0}=0$ in our case).
	\end{proof}
	
	\begin{rem}
		The definition of the polynomials $A_{n}^{s}(x)$ is equivalent to 
		\begin{align*}
			A_{n}^{s}(x)=  \sum_{k=1}^{n} (n+k)^{s} \left\lfloor\!\!\left\lfloor\!
			\begin{matrix}
				n \\ k
			\end{matrix}
			\!\right\rfloor\!\!\right\rfloor_{s} x^{n+k},
		\end{align*} since the term corresponding to $k=0$ equals $0$.  
	\end{rem}
	
	\begin{lem}
		For all $n\geq1$, the polynomials $A_{n}^{s}(x)$ satisfy
		\begin{align*}
			A_{n}^{s} = \sum_{i=0}^{s} \sttwo{s}{i} x^{i} \frac{d^{i}}{dx^{i}} Q_{n}^{s}(x).
		\end{align*}
	\end{lem}
	
	\begin{proof}
		Since $\sttwo{n}{k} = \left\{\mkern-4mu\braceVectorstack{n\\k}\mkern-4mu\right\}_{1}$, the Stirling numbers of the second kind satisfy
		\begin{align}
			\label{stirpolid}
			x^{n} = \sum_{k=0}^{n} \sttwo{n}{k} x^{\underline{k}}.
		\end{align}
		It follows from (\ref{stirpolid}) that
		\begin{align}
			A_{n}^{s}(x) &= \sum_{k=0}^{n} (n+k)^{s} \left\lfloor\!\!\left\lfloor\!
			\begin{matrix}
				n \\ k
			\end{matrix}
			\!\right\rfloor\!\!\right\rfloor_{s} x^{n+k}  \nonumber \\
			&= \sum_{k=0}^{n} \left\lfloor\!\!\left\lfloor\!
			\begin{matrix}
				n \\ k
			\end{matrix}
			\!\right\rfloor\!\!\right\rfloor_{s} (n+k)^{s} x^{n+k}  \nonumber \\
			&= \sum_{k=0}^{n} \left\lfloor\!\!\left\lfloor\!
			\begin{matrix}
				n \\ k
			\end{matrix}
			\!\right\rfloor\!\!\right\rfloor_{s} \sum_{i=0}^{s} \sttwo{s}{i} (n+k)^{\underline{i}} x^{n+k}  \nonumber \\
			&= \sum_{i=0}^{s} \sttwo{s}{i} x^{i} \sum_{k=0}^{n} (n+k)^{\underline{i}} \left\lfloor\!\!\left\lfloor\!
			\begin{matrix}
				n \\ k
			\end{matrix}
			\!\right\rfloor\!\!\right\rfloor_{s} x^{n+k-i}  \nonumber \\
			&= \sum_{i=0}^{s} \sttwo{s}{i} x^{i} \frac{d^{i}}{dx^{i}} Q_{n}^{s}(x). 	\label{ansx}
		\end{align}
	\end{proof} Therefore, (\ref{ddQ}) can be written as
	\begin{align*}
		Q_{n+1}^{s}(x) &= x^{2} Q_{n}^{s}(x) + xA_{n}^{s}(x) \\
		&= x^{2} Q_{n}^{s}(x) + x \left(\sum_{i=0}^{s} \sttwo{s}{i} x^{i} \frac{d^{i}}{dx^{i}} Q_{n}^{s}(x) \right) \\
		&=  x^{2} Q_{n}^{s}(x) + \sum_{i=0}^{s} \sttwo{s}{i} x^{i+1} \frac{d^{i}}{dx^{i}} Q_{n}^{s}(x),
	\end{align*} which is equivalent to
	
	\begin{equation}
		\label{recQ}
		Q_{n+1}^{s}(x) = x^{2} Q_{n}^{s}(x) + \sum_{i=1}^{s} \sttwo{s}{i} x^{i+1} \frac{d^{i}}{dx^{i}} Q_{n}^{s}(x),
	\end{equation} since $s\geq 1$ by the definition of $Q_{n}^{s}(x)$ and $A_{n}^{s}(x)$ (the term corresponding to $i=0$ in (\ref{ansx}) is nonzero only when $s=i=0$).
	
	\begin{cor}
		In terms of the Lah polynomials with higher level, (\ref{recQ}) can be written as
		\begin{equation}
			\label{row2}
			L_{n+1}^{s}(x) = xL_{n}^{s}(x) + \sum_{i=1}^{s} \sttwo{s}{i} x^{i-n} \frac{d^{i}}{dx^{i}} \left(x^{n} L_{n}^{s} (x) \right).
		\end{equation}
	\end{cor}
	
	\begin{proof}
		In equation (\ref{recQ}), we consider that $x^{n+1} L_{n+1}^{s}(x) = Q_{n+1}^{s} (x)$ and we get 
		\begin{align*}
			x^{n+1} L_{n+1}^{s} (x) &= x^{2} x^{n} L_{n}^{s} (x) + \sum_{i=1}^{s} \sttwo{s}{i} x^{i+1} \frac{d^i}{dx^i} \left( x^{n} L_{n}^{s} (x) \right) \\
			&= x^{n+2} L_{n}^{s} (x) + \sum_{i=1}^{s} \sttwo{s}{i} x^{i+1} \frac{d^i}{dx^i} \left( x^{n} L_{n}^{s} (x) \right).
		\end{align*}
		By dividing the equation by $x^{n+1}$, we get (\ref{row2}).
	\end{proof}

	\section{The Lah Numbers of Order $s$}
	
	In this section, we generalize Luschny's Lah numbers of order 2 and 3 \cite[A268434, A269946]{SloaneOEIS} to the Lah numbers of order $s$ with the relation involving Stirling numbers with higher level of both kinds, the recurrence  relation and the polynomial identity including falling and rising factorials with higher level. We also study the Lah polynomials of order $s$. 
	
	\subsection{Definition, Recurrence Relation and Some Special Cases}
	
	\begin{defn}
		\label{defslah}
		The \textit{Lah numbers of order $s$},  $\left\lceil\!\!\left\lceil\!
		\begin{matrix}
			n \\ k
		\end{matrix}
		\!\right\rceil\!\!\right\rceil_{s}$ are for $n, k \in \mathbb{N}_{0}$ and $s \in \mathbb{N}$ defined by the relation
		\begin{equation}
			\label{slahdef}
			\left\lceil\!\!\left\lceil\!
			\begin{matrix}
				n \\ k
			\end{matrix}
			\!\right\rceil\!\!\right\rceil_{s} = \sum_{j=k}^{n}  \begin{NiceArrayWithDelims}{\llbracket}{\rrbracket}{c}
				n \\ j
			\end{NiceArrayWithDelims}_{s} \left\{\mkern-4mu\braceVectorstack{j\\k}\mkern-4mu\right\}_{s}
		\end{equation}
		with $\left\lceil\!\!\left\lceil\!
		\begin{matrix}
			n \\ k
		\end{matrix}
		\!\right\rceil\!\!\right\rceil_{s}=0$ for $n<k$.
	\end{defn}
	
	\begin{ex}
		Let us give an example for $n=4$, $k=3$ and $s=2$. Proceeding with the previous definition, we get 
		\begin{align*}
			\left\lceil\!\!\left\lceil\!
			\begin{matrix}
				4 \\ 3
			\end{matrix}
			\!\right\rceil\!\!\right\rceil_{2}  & \left.= \right.  \sum_{j=3}^{4} \begin{NiceArrayWithDelims}{\llbracket}{\rrbracket}{c}
				4 \\ j
			\end{NiceArrayWithDelims}_{2} \left\{\mkern-4mu\braceVectorstack{j\\3}\mkern-4mu\right\}_{2} \\
			&=  \begin{NiceArrayWithDelims}{\llbracket}{\rrbracket}{c}
				4 \\ 3
			\end{NiceArrayWithDelims}_{2} \left\{\mkern-4mu\braceVectorstack{3\\3}\mkern-4mu\right\}_{2} + \begin{NiceArrayWithDelims}{\llbracket}{\rrbracket}{c}
				4 \\ 4
			\end{NiceArrayWithDelims}_{2} \left\{\mkern-4mu\braceVectorstack{4\\3}\mkern-4mu\right\}_{2} \\
			&= \begin{NiceArrayWithDelims}{\llbracket}{\rrbracket}{c}
				4 \\ 3
			\end{NiceArrayWithDelims}_{2} + \left\{\mkern-4mu\braceVectorstack{4\\3}\mkern-4mu\right\}_{2} = 28.
		\end{align*}
	\end{ex}

	\begin{thm}
		For $n, k \geq 1$ and $s \in \mathbb{N}$, the Lah numbers of order $s$ satisfy the recurrence 
		\begin{equation}
			\label{slahrec}
			\left\lceil\!\!\left\lceil\!
			\begin{matrix}
				n \\ k
			\end{matrix}
			\!\right\rceil\!\!\right\rceil_{s} = \left\lceil\!\!\left\lceil\!
			\begin{matrix}
				n-1 \\ k-1
			\end{matrix}
			\!\right\rceil\!\!\right\rceil_{s} + ((n-1)^{s} + k^{s}) \left\lceil\!\!\left\lceil\!
			\begin{matrix}
				n-1 \\ k
			\end{matrix}
			\!\right\rceil\!\!\right\rceil_{s}
		\end{equation} with boundary conditions
		
		\begin{align*}
			\left\lceil\!\!\left\lceil\!
			\begin{matrix}
				n \\ n
			\end{matrix}
			\!\right\rceil\!\!\right\rceil_{s} = 1 \hspace{2mm} \text{for $n \geq 0$,}
		\end{align*}
		and
		\begin{align*}
			\left\lceil\!\!\left\lceil\!
			\begin{matrix}
				n \\ 0
			\end{matrix}
			\!\right\rceil\!\!\right\rceil_{s} = \left\lceil\!\!\left\lceil\!
			\begin{matrix}
				0 \\ n
			\end{matrix}
			\!\right\rceil\!\!\right\rceil_{s} =  0 \hspace{2mm} \text{for $n>0$.}
		\end{align*}
	\end{thm}
	
	\begin{proof}
		We need to show that the right-hand side of (\ref{slahdef}) satisfies the recurrence relation (\ref{slahrec}). We proceed by applying recurrence relations (\ref{sstirling1rec}) and (\ref{sstirling2rec}). 
		
		\begin{align*}
			\left\lceil\!\!\left\lceil\!
			\begin{matrix}
				n \\ k
			\end{matrix}
			\!\right\rceil\!\!\right\rceil_{s} &= \sum_{j=k}^{n}  \begin{NiceArrayWithDelims}{\llbracket}{\rrbracket}{c}
				n \\ j
			\end{NiceArrayWithDelims}_{s} \left\{\mkern-4mu\braceVectorstack{j\\k}\mkern-4mu\right\}_{s} \\
			&= \sum_{j=k}^{n}  \left( \begin{NiceArrayWithDelims}{\llbracket}{\rrbracket}{c}
				n-1 \\ j-1
			\end{NiceArrayWithDelims}_{s} + (n-1)^{s} \begin{NiceArrayWithDelims}{\llbracket}{\rrbracket}{c}
				n-1 \\ j
			\end{NiceArrayWithDelims}_{s} \right) \left\{\mkern-4mu\braceVectorstack{j\\k}\mkern-4mu\right\}_{s} \\
			&= \sum_{j=k}^{n} \begin{NiceArrayWithDelims}{\llbracket}{\rrbracket}{c}
				n-1 \\ j-1
			\end{NiceArrayWithDelims}_{s} \left\{\mkern-4mu\braceVectorstack{j\\k}\mkern-4mu\right\}_{s} + (n-1)^{s} \sum_{j=k}^{n-1}\begin{NiceArrayWithDelims}{\llbracket}{\rrbracket}{c}
				n-1 \\ j
			\end{NiceArrayWithDelims}_{s} \left\{\mkern-4mu\braceVectorstack{j\\k}\mkern-4mu\right\}_{s} + (n-1)^{s} \begin{NiceArrayWithDelims}{\llbracket}{\rrbracket}{c}
				n-1 \\ n
			\end{NiceArrayWithDelims}_{s} \left\{\mkern-4mu\braceVectorstack{n\\k}\mkern-4mu\right\}_{s} \\
			&= \sum_{j=k}^{n} \begin{NiceArrayWithDelims}{\llbracket}{\rrbracket}{c}
				n-1 \\ j-1
			\end{NiceArrayWithDelims}_{s} \left( \left\{\mkern-4mu\braceVectorstack{j-1\\k-1}\mkern-4mu\right\}_{s} + k^{s} \left\{\mkern-4mu\braceVectorstack{j-1\\k}\mkern-4mu\right\}_{s} \right) + (n-1)^{s} \left\lceil\!\!\left\lceil\!
			\begin{matrix}
				n-1 \\ k
			\end{matrix}
			\!\right\rceil\!\!\right\rceil_{s} \\
			&= \sum_{j=k}^{n} \begin{NiceArrayWithDelims}{\llbracket}{\rrbracket}{c}
				n-1 \\ j-1
			\end{NiceArrayWithDelims}_{s} \left\{\mkern-4mu\braceVectorstack{j-1\\k-1}\mkern-4mu\right\}_{s} + k^{s} \sum_{j=k}^{n} \begin{NiceArrayWithDelims}{\llbracket}{\rrbracket}{c}
				n-1 \\ j-1
			\end{NiceArrayWithDelims}_{s} \left\{\mkern-4mu\braceVectorstack{j-1\\k}\mkern-4mu\right\}_{s}  + (n-1)^{s} \left\lceil\!\!\left\lceil\!
			\begin{matrix}
				n-1 \\ k
			\end{matrix}
			\!\right\rceil\!\!\right\rceil_{s} \\
			&= \sum_{j=k-1}^{n-1} \begin{NiceArrayWithDelims}{\llbracket}{\rrbracket}{c}
				n-1 \\ j
			\end{NiceArrayWithDelims}_{s} \left\{\mkern-4mu\braceVectorstack{j\\k-1}\mkern-4mu\right\}_{s}  + k^{s} \sum_{j=k-1}^{n-1} \begin{NiceArrayWithDelims}{\llbracket}{\rrbracket}{c}
				n-1 \\ j
			\end{NiceArrayWithDelims}_{s} \left\{\mkern-4mu\braceVectorstack{j\\k}\mkern-4mu\right\}_{s}   + (n-1)^{s} \left\lceil\!\!\left\lceil\!
			\begin{matrix}
				n-1 \\ k
			\end{matrix}
			\!\right\rceil\!\!\right\rceil_{s} \\
			&= \left\lceil\!\!\left\lceil\!
			\begin{matrix}
				n-1 \\ k-1
			\end{matrix}
			\!\right\rceil\!\!\right\rceil_{s} + k^{s} \begin{NiceArrayWithDelims}{\llbracket}{\rrbracket}{c}
				n-1 \\ k-1
			\end{NiceArrayWithDelims}_{s} \left\{\mkern-4mu\braceVectorstack{k-1\\k}\mkern-4mu\right\}_{s} + k^{s} \sum_{j=k}^{n-1} \begin{NiceArrayWithDelims}{\llbracket}{\rrbracket}{c}
				n-1 \\ j
			\end{NiceArrayWithDelims}_{s} \left\{\mkern-4mu\braceVectorstack{j\\k}\mkern-4mu\right\}_{s} + (n-1)^{s} \left\lceil\!\!\left\lceil\!
			\begin{matrix}
				n-1 \\ k
			\end{matrix}
			\!\right\rceil\!\!\right\rceil_{s} \\
			&= \left\lceil\!\!\left\lceil\!
			\begin{matrix}
				n-1 \\ k-1
			\end{matrix}
			\!\right\rceil\!\!\right\rceil_{s} + k^{s}\left\lceil\!\!\left\lceil\!
			\begin{matrix}
				n-1 \\ k
			\end{matrix}
			\!\right\rceil\!\!\right\rceil_{s} + (n-1)^{s} \left\lceil\!\!\left\lceil\!
			\begin{matrix}
				n-1 \\ k
			\end{matrix}
			\!\right\rceil\!\!\right\rceil_{s}    \\
			&= \left\lceil\!\!\left\lceil\!
			\begin{matrix}
				n-1 \\ k-1
			\end{matrix}
			\!\right\rceil\!\!\right\rceil_{s} + ((n-1)^{s} + k^{s}) \left\lceil\!\!\left\lceil\!
			\begin{matrix}
				n-1 \\ k
			\end{matrix}
			\!\right\rceil\!\!\right\rceil_{s}.
		\end{align*} Now, we show that the boundary conditions also match. Applying (\ref{cases1a}) and (\ref{cases1b}), we get
		
		\begin{align*}
			&\left\lceil\!\!\left\lceil\!
			\begin{matrix}
				n \\ n
			\end{matrix}
			\!\right\rceil\!\!\right\rceil_{s} = \sum_{j=n}^{n} \begin{NiceArrayWithDelims}{\llbracket}{\rrbracket}{c}
				n \\ j
			\end{NiceArrayWithDelims}_{s} \left\{\mkern-4mu\braceVectorstack{j\\n}\mkern-4mu\right\}_{s} = \begin{NiceArrayWithDelims}{\llbracket}{\rrbracket}{c}
				n \\ n
			\end{NiceArrayWithDelims}_{s} \left\{\mkern-4mu\braceVectorstack{n\\n}\mkern-4mu\right\}_{s} = 1 \hspace{2mm} \text{for $n \geq 0$,} \\
			&
			\left\lceil\!\!\left\lceil\!
			\begin{matrix}
				n \\ 0
			\end{matrix}
			\!\right\rceil\!\!\right\rceil_{s} = \sum_{j=0}^{n} \begin{NiceArrayWithDelims}{\llbracket}{\rrbracket}{c}
				n \\ j
			\end{NiceArrayWithDelims}_{s} \left\{\mkern-4mu\braceVectorstack{j\\0}\mkern-4mu\right\}_{s} = \begin{NiceArrayWithDelims}{\llbracket}{\rrbracket}{c}
				n \\ 0
			\end{NiceArrayWithDelims}_{s} \left\{\mkern-4mu\braceVectorstack{0\\0}\mkern-4mu\right\}_{s} = \begin{NiceArrayWithDelims}{\llbracket}{\rrbracket}{c}
				n \\ 0
			\end{NiceArrayWithDelims}_{s} = 0  \hspace{2mm} \text{for $n$ $>$ 0.} 
		\end{align*} That is, the recurrence (\ref{slahrec}) holds for all $n, k \geq 1$ and $s \in \mathbb{N}$.
	\end{proof}
	
	\begin{thm}
		(Some special cases)
		\begin{align}
			&\label{special1} \left\lceil\!\!\left\lceil\!
			\begin{matrix}
				2 \\ 1
			\end{matrix}
			\!\right\rceil\!\!\right\rceil_{s} = 2 \\
			& \label{special2}\left\lceil\!\!\left\lceil\!
			\begin{matrix}
				3 \\ 1
			\end{matrix}
			\!\right\rceil\!\!\right\rceil_{s} = \left\lceil\!\!\left\lceil\!
			\begin{matrix}
				3 \\ 2
			\end{matrix}
			\!\right\rceil\!\!\right\rceil_{s} = 2^{s+1} +2 \\
			& \label{k=n-1(2)}\left\lceil\!\!\left\lceil\!
			\begin{matrix}
				n \\ n-1
			\end{matrix}
			\!\right\rceil\!\!\right\rceil_{s} = \sum_{j=1}^{n-1} 2j^{s}    \\
			& \label{k=1} \left\lceil\!\!\left\lceil\!
			\begin{matrix}
				n \\ 1
			\end{matrix}
			\!\right\rceil\!\!\right\rceil_{s} = \prod_{i=1}^{n} ((n-i)^{s}+1).
		\end{align}
	\end{thm}

	\begin{proof}
		First, we prove (\ref{special1}). Using (\ref{slahrec}), we get
		\begin{align*}
			\left\lceil\!\!\left\lceil\!
			\begin{matrix}
				2 \\ 1
			\end{matrix}
			\!\right\rceil\!\!\right\rceil_{s} = 	\left\lceil\!\!\left\lceil\!
			\begin{matrix}
				1\\ 0
			\end{matrix} 
			\!\right\rceil\!\!\right\rceil_{s} + (1^\textbf{s} + 1^{s}) 	\left\lceil\!\!\left\lceil\!
			\begin{matrix}
				1 \\ 1
			\end{matrix}
			\!\right\rceil\!\!\right\rceil_{s} = 1^{s} + 1^{s} = 2.
		\end{align*} \\
		Next, we prove (\ref{special2}). Using (\ref{slahrec}), we get
		\begin{align*}
			\left\lceil\!\!\left\lceil\!
			\begin{matrix}
				3 \\ 1
			\end{matrix}
			\!\right\rceil\!\!\right\rceil_{s} = 	\left\lceil\!\!\left\lceil\!
			\begin{matrix}
				2\\ 0
			\end{matrix} 
			\!\right\rceil\!\!\right\rceil_{s} + (2^\textbf{s} + 1^{s}) 	\left\lceil\!\!\left\lceil\!
			\begin{matrix}
				2 \\ 1
			\end{matrix}
			\!\right\rceil\!\!\right\rceil_{s} = 2 (2^{s} + 1^{s}) =  2^{s+1} +2,
		\end{align*}
		and
		\begin{align*}
			\left\lceil\!\!\left\lceil\!
			\begin{matrix}
				3 \\ 2
			\end{matrix}
			\!\right\rceil\!\!\right\rceil_{s} = 	\left\lceil\!\!\left\lceil\!
			\begin{matrix}
				2\\ 1
			\end{matrix} 
			\!\right\rceil\!\!\right\rceil_{s} + (2^\textbf{s} + 2^{s}) 	\left\lceil\!\!\left\lceil\!
			\begin{matrix}
				2 \\ 2
			\end{matrix}
			\!\right\rceil\!\!\right\rceil_{s} = 2 + (2^{s} + 2^{s}) =  2\cdot 2^{s} + 2= 2^{s+1} + 2.
		\end{align*} \\
		Now, we prove (\ref{k=n-1(2)}). By (\ref{slahrec}) with $k=n-1$ we have
		\begin{align}
			\left\lceil\!\!\left\lceil\!
			\begin{matrix}
				n \\ n-1
			\end{matrix}
			\!\right\rceil\!\!\right\rceil_{s}\ &=\ \left\lceil\!\!\left\lceil\!
			\begin{matrix}
				n-1 \\ n-2
			\end{matrix}
			\!\right\rceil\!\!\right\rceil_{s} + ((n-1)^{s} + (n-1)^{s}) \left\lceil\!\!\left\lceil\!
			\begin{matrix}
				n-1 \\ n-1
			\end{matrix}
			\!\right\rceil\!\!\right\rceil_{s} \nonumber \\
			&= \left\lceil\!\!\left\lceil\!
			\begin{matrix}
				n-1 \\ n-2
			\end{matrix}
			\!\right\rceil\!\!\right\rceil_{s} + 2 (n-1)^{s} 	\label{slahboundaryproof}
		\end{align}
		since $\left\lceil\!\!\left\lceil\!
		\begin{matrix}
			n-1 \\ n-1
		\end{matrix}
		\!\right\rceil\!\!\right\rceil_{s}=1$. Using (\ref{slahboundaryproof}) repeatedly, we obtain  
		\begin{align*}
			\left\lceil\!\!\left\lceil\!
			\begin{matrix}
				n \\ n-1
			\end{matrix}
			\!\right\rceil\!\!\right\rceil_{s}\ &=\ \underline{\left\lceil\!\!\left\lceil\!
				\begin{matrix}
					n-1 \\ n-2
				\end{matrix}
				\!\right\rceil\!\!\right\rceil_{s}} +\ 2(n-1)^{s} \\
			\underline{\left\lceil\!\!\left\lceil\!
				\begin{matrix}
					n-1 \\ n-2
				\end{matrix}
				\!\right\rceil\!\!\right\rceil_{s}}\ &=\ \underline{\left\lceil\!\!\left\lceil\!
				\begin{matrix}
					n-2 \\ n-3
				\end{matrix}
				\!\right\rceil\!\!\right\rceil_{s}} +\ 2(n-2)^{s} \\
			&\ \,\vdots \\
			\underline{\left\lceil\!\!\left\lceil\!
				\begin{matrix}
					2 \\ 1
				\end{matrix}
				\!\right\rceil\!\!\right\rceil_{s}}\ &=\ \left\lceil\!\!\left\lceil\!
			\begin{matrix}
				1 \\ 0
			\end{matrix}
			\!\right\rceil\!\!\right\rceil_{s} +\ 2 \cdot 1.
		\end{align*}
		Summing these equations and cancelling all the underlined terms, we find
		\begin{align*}
			\left\lceil\!\!\left\lceil\!
			\begin{matrix}
				n \\ n-1
			\end{matrix}
			\!\right\rceil\!\!\right\rceil_{s} = \left\lceil\!\!\left\lceil\!
			\begin{matrix}
				1 \\ 0
			\end{matrix}
			\!\right\rceil\!\!\right\rceil_{s} + \sum_{j=1}^{n-1} 2(n-j)^{s} = \sum_{j=1}^{n-1} 2(n-j)^{s} = \sum_{j=1}^{n-1} 2j^{s}.
		\end{align*} \\ 
		Finally, we prove (\ref{k=1}). By (\ref{slahrec}) with $k=1$ we have
		\begin{align}
			\left\lceil\!\!\left\lceil\!
			\begin{matrix}
				n \\ 1
			\end{matrix}
			\!\right\rceil\!\!\right\rceil_{s}\ &=\ \left\lceil\!\!\left\lceil\!
			\begin{matrix}
				n-1 \\ 0
			\end{matrix}
			\!\right\rceil\!\!\right\rceil_{s} + ((n-1)^{s} + 1^{s}) \left\lceil\!\!\left\lceil\!
			\begin{matrix}
				n-1 \\ 1
			\end{matrix}
			\!\right\rceil\!\!\right\rceil_{s} \nonumber \\
			&= ((n-1)^{s} +1) \left\lceil\!\!\left\lceil\!
			\begin{matrix}
				n-1 \\ 1
			\end{matrix}
			\!\right\rceil\!\!\right\rceil_{s} 	\label{slahboundary2proof}
		\end{align}
		since $\left\lceil\!\!\left\lceil\!
		\begin{matrix}
			n-1 \\ 0
		\end{matrix}
		\!\right\rceil\!\!\right\rceil_{s} = 0$. Using (\ref{slahboundary2proof}) repeatedly, we obtain
		\begin{align*}
			\left\lceil\!\!\left\lceil\!
			\begin{matrix}
				n \\ 1
			\end{matrix}
			\!\right\rceil\!\!\right\rceil_{s}\ &=\ ((n-1)^{s}+1) 	\left\lceil\!\!\left\lceil\!
			\begin{matrix}
				n-1 \\ 1
			\end{matrix}
			\!\right\rceil\!\!\right\rceil_{s} \\
			\left\lceil\!\!\left\lceil\!
			\begin{matrix}
				n-1 \\ 1
			\end{matrix}
			\!\right\rceil\!\!\right\rceil_{s}\ &=\ ((n-2)^{s}+1) 	\left\lceil\!\!\left\lceil\!
			\begin{matrix}
				n-2 \\ 1
			\end{matrix}
			\!\right\rceil\!\!\right\rceil_{s} \\
			&\ \,\vdots \\
			\left\lceil\!\!\left\lceil\!
			\begin{matrix}
				1 \\ 1
			\end{matrix}
			\!\right\rceil\!\!\right\rceil_{s}\ &=\ 1.
		\end{align*}
		By writing all terms from $\left\lceil\!\!\left\lceil\!
		\begin{matrix}
			n \\ 1
		\end{matrix}
		\!\right\rceil\!\!\right\rceil_{s}$ to $\left\lceil\!\!\left\lceil\!
		\begin{matrix}
			1 \\ 1
		\end{matrix}
		\!\right\rceil\!\!\right\rceil_{s}$ in this way, we obtain
		\begin{align*}
			\left\lceil\!\!\left\lceil\!
			\begin{matrix}
				n \\ n-1
			\end{matrix}
			\!\right\rceil\!\!\right\rceil_{s} = ((n-1)^{s}+1) ((n-2)^{s}+1) \cdots ((n-n)^{s}+1) =  \prod_{i=1}^{n} ((n-i)^{s}+1).
		\end{align*}
	\end{proof}
	
	\begin{rem}
		Note that for $s=1$, we get the ordinary Lah numbers.
	\end{rem}
	
	Using the recurrence relation, we get tables of values (see Tables \ref{table1s} - \ref{table3s}). For example, $\left\lceil\!\!\left\lceil\!
	\begin{matrix}
		6\\ 2
	\end{matrix}
	\!\right\rceil\!\!\right\rceil_{25}$ is 
	
	\begin{align*}
		19080000046878387911728136488155674014369724428110000.
	\end{align*}	
	\begin{table}[ht]
		\begin{center}
			\begin{tabular}{|c|c|c|c|c|c|c|c|}
				\hline
				\diagbox[width=2em]{$n$}{$k$} & 0& 1 &2& 3 & 4 & 5 &6  \\
				\hline
				0 & 1 & 0 & 0 & 0 &  0 & 0 & 0  \\
				\hline
				1 & 0 & 1 & 0 & 0 & 0 & 0 & 0  \\
				\hline
				2 & 0 & 2 & 1 & 0 & 0 &  0 & 0   \\
				\hline
				3 & 0 & 10 & 10 & 1 & 0 & 0 & 0  \\
				\hline
				4 & 0 & 100 & 140 & 28 & 1 & 0 & 0  \\
				\hline
				5 & 0 & 1700 & 2900 & 840 & 60 & 1 & 0  \\
				\hline
				6 & 0 & 44200 & 85800 & 31460 & 3300 & 110 & 1  \\
				\hline
			\end{tabular}
		\end{center}
		\caption{The Lah numbers of order $2$, $\left\lceil\!\!\left\lceil\!
			\begin{matrix}
				n \\ k
			\end{matrix}
			\!\right\rceil\!\!\right\rceil_{2}$ for $0 \le n, k\le 6$.}
		\label{table1s}
	\end{table}

	\begin{table}[ht]
		\begin{center}
			\begin{tabular}{|c|c|c|c|c|c|c|c|}
				\hline
				\diagbox[width=2em]{$n$}{$k$} & 0& 1 &2& 3 & 4 & 5 &6  \\
				\hline
				0 & 1 & 0 & 0 & 0 &  0 & 0 & 0  \\
				\hline
				1 & 0 & 1 & 0 & 0 & 0 & 0 & 0  \\
				\hline
				2 & 0 & 2 & 1 & 0 & 0 &  0 & 0   \\
				\hline
				3 & 0 & 18 & 18 & 1 & 0 & 0 & 0  \\
				\hline
				4 & 0 & 504 & 648 & 72 & 1 & 0 & 0  \\
				\hline
				5 & 0 & 32760 & 47160 & 7200 & 200 & 1 & 0  \\
				\hline
				6 & 0 & 4127760 & 6305040 & 1141560 & 45000 & 450 & 1  \\
				\hline
			\end{tabular}
		\end{center}
		\caption{The Lah numbers of order $3$, $\left\lceil\!\!\left\lceil\!
			\begin{matrix}
				n \\ k
			\end{matrix}
			\!\right\rceil\!\!\right\rceil_{3}$ for $0 \le n, k\le 6$.}
		\label{table2s}
	\end{table}

	\begin{table}[ht]
		\begin{center}
			\begin{tabular}{|c|c|c|c|c|c|c|c|}
				\hline
				\diagbox[width=2em]{$n$}{$k$} & 0& 1 &2& 3 & 4 & 5 &6  \\
				\hline
				0 & 1 & 0 & 0 & 0 &  0 & 0 & 0  \\
				\hline
				1 & 0 & 1 & 0 & 0 & 0 & 0 & 0  \\
				\hline
				2 & 0 & 2 & 1 & 0 & 0 &  0 & 0   \\
				\hline
				3 & 0 & 34 & 34 & 1 & 0 & 0 & 0  \\
				\hline
				4 & 0 & 2788 & 3332 & 196 & 1 & 0 & 0  \\
				\hline
				5 & 0 & 716516 & 909092 & 69384 & 708 & 1 & 0  \\
				\hline
				6 & 0 & 448539016 & 583444488 & 49894196 & 693132 & 1958 & 1  \\
				\hline
			\end{tabular}
		\end{center}
		\caption{The Lah numbers of order $4$, $\left\lceil\!\!\left\lceil\!
			\begin{matrix}
				n \\ k
			\end{matrix}
			\!\right\rceil\!\!\right\rceil_{4}$ for $0 \le n, k\le 6$.}
		\label{table3s}
	\end{table}

	\subsection{Relations via Falling and Rising Factorials with Higher Level}	
	
	\begin{thm}
		\label{risingidentity}
		For $n \geq 0$, the Lah numbers of order $s$ satisfy the polynomial identity 
		\begin{equation}
			\label{rising}
			x^{\overline{\overline{n}}}_{s}= \sum_{k=0}^{n} \left\lceil\!\!\left\lceil\!
			\begin{matrix}
				n \\ k
			\end{matrix}
			\!\right\rceil\!\!\right\rceil_{s} x^{\underline{\underline{k}}}_{s}.
		\end{equation}
	\end{thm}
	
	\begin{proof}
		We proceed by induction on $n$. The equality clearly holds for $n=0$. We now assume that the identity holds for some $n \geq 0$ and prove it for $n+1$. Using the recurrence relation (\ref{slahrec}) and induction hypothesis, we thus obtain
		\begin{align*}
			x^{\overline{\overline{n+1}}}_{s} &= (x+n^{s}) x^{\overline{\overline{n}}}_{s} \\
			&= ((x-k^{s})+(n^{s}+k^{s})) x^{\overline{\overline{n}}}_{s} \\
			&=  ((x-k^{s})+(n^{s}+k^{s})) \left( \sum_{k=0}^{n} \left\lceil\!\!\left\lceil\!
			\begin{matrix}
				n \\ k
			\end{matrix}
			\!\right\rceil\!\!\right\rceil_{s} x^{\underline{\underline{k}}}_{s} \right) \\
			&= \sum_{k=0}^{n} \left( (x-k^{s}) \left\lceil\!\!\left\lceil\!
			\begin{matrix}
				n \\ k
			\end{matrix}
			\!\right\rceil\!\!\right\rceil_{s} x^{\underline{\underline{k}}}_{s} + (n^{s}+k^{s}) \left\lceil\!\!\left\lceil\!
			\begin{matrix}
				n \\ k
			\end{matrix}
			\!\right\rceil\!\!\right\rceil_{s} x^{\underline{\underline{k}}}_{s} \right) \\
			&= \sum_{k=0}^{n} (x-k^{s}) \left\lceil\!\!\left\lceil\!
			\begin{matrix}
				n \\ k
			\end{matrix}
			\!\right\rceil\!\!\right\rceil_{s} x^{\underline{\underline{k}}}_{s} + \sum_{k=0}^{n} (n^{s} + k^{s}) \left\lceil\!\!\left\lceil\!
			\begin{matrix}
				n \\ k
			\end{matrix}
			\!\right\rceil\!\!\right\rceil_{s} x^{\underline{\underline{k}}}_{s} \\
			&= \sum_{k=0}^{n} \left\lceil\!\!\left\lceil\!
			\begin{matrix}
				n \\ k
			\end{matrix}
			\!\right\rceil\!\!\right\rceil_{s} x^{\underline{\underline{k+1}}}_{s} +  \sum_{k=0}^{n} (n^{s} + k^{s})  \left\lceil\!\!\left\lceil\!
			\begin{matrix}
				n \\ k
			\end{matrix}
			\!\right\rceil\!\!\right\rceil_{s} x^{\underline{\underline{k}}}_{s}  \\
			&= \sum_{k=1}^{n+1} \left\lceil\!\!\left\lceil\!
			\begin{matrix}
				n \\ k-1
			\end{matrix}
			\!\right\rceil\!\!\right\rceil_{s} x^{\underline{\underline{k}}}_{s} +  \sum_{k=0}^{n} (n^{s} + k^{s})  \left\lceil\!\!\left\lceil\!
			\begin{matrix}
				n \\ k
			\end{matrix}
			\!\right\rceil\!\!\right\rceil_{s} x^{\underline{\underline{k}}}_{s}  \\
			&= \sum_{k=1}^{n+1} \left( \left\lceil\!\!\left\lceil\!
			\begin{matrix}
				n \\ k-1
			\end{matrix}
			\!\right\rceil\!\!\right\rceil_{s} + (n^{s} + k^{s}) \left\lceil\!\!\left\lceil\!
			\begin{matrix}
				n \\ k
			\end{matrix}
			\!\right\rceil\!\!\right\rceil_{s} \right) x^{\underline{\underline{k}}}_{s} \\
			&= \sum_{k=1}^{n+1} \left\lceil\!\!\left\lceil\!
			\begin{matrix}
				n+1 \\ k
			\end{matrix}
			\!\right\rceil\!\!\right\rceil_{s} x^{\underline{\underline{k}}}_{s} = \sum_{k=0}^{n+1} \left\lceil\!\!\left\lceil\!
			\begin{matrix}
				n+1 \\ k
			\end{matrix}
			\!\right\rceil\!\!\right\rceil_{s} x^{\underline{\underline{k}}}_{s},
		\end{align*}
		which completes the induction step. Thus, by induction, the identity holds for all $n \geq 0$.
	\end{proof}
	
	\begin{rem}
		The polynomial identity (\ref{rising}) can be seen as a row generating function for the Lah numbers of order $s$. 
	\end{rem}
	
	\begin{ex}
		Let us give an example of Theorem \ref{risingidentity} for $n=4$ and $s=2$.
		\begin{align*}
			x^{\overline{\overline{4}}}_{2} &\left.= \right.  \sum_{k=0}^{4} \left\lceil\!\!\left\lceil\!
			\begin{matrix}
				4 \\ k
			\end{matrix}
			\!\right\rceil\!\!\right\rceil_{2} x^{\underline{\underline{k}}}_{2} \\
			&=\left\lceil\!\!\left\lceil\!
			\begin{matrix}
				4 \\ 0
			\end{matrix}
			\!\right\rceil\!\!\right\rceil_{2} x^{\underline{\underline{0}}}_{2} +  \left\lceil\!\!\left\lceil\!
			\begin{matrix}
				4 \\ 1
			\end{matrix}
			\!\right\rceil\!\!\right\rceil_{2} x^{\underline{\underline{1}}}_{2} + \left\lceil\!\!\left\lceil\!
			\begin{matrix}
				4 \\ 2
			\end{matrix}
			\!\right\rceil\!\!\right\rceil_{2} x^{\underline{\underline{2}}}_{2} + \left\lceil\!\!\left\lceil\!
			\begin{matrix}
				4 \\ 3
			\end{matrix}
			\!\right\rceil\!\!\right\rceil_{2} x^{\underline{\underline{3}}}_{2} + \left\lceil\!\!\left\lceil\!
			\begin{matrix}
				4 \\ 4
			\end{matrix}
			\!\right\rceil\!\!\right\rceil_{2} x^{\underline{\underline{4}}}_{2}.
		\end{align*} 
		We compute the values of the Lah numbers of order $s$ using the recurrence relation (\ref{slahrec}): $\left\lceil\!\!\left\lceil\!
		\begin{matrix}
			4 \\ 0 
		\end{matrix}
		\!\right\rceil\!\!\right\rceil_{2}= 0, \left\lceil\!\!\left\lceil\!
		\begin{matrix}
			4 \\ 1
		\end{matrix}
		\!\right\rceil\!\!\right\rceil_{2} = 100, \left\lceil\!\!\left\lceil\!
		\begin{matrix}
			4 \\ 2
		\end{matrix}
		\!\right\rceil\!\!\right\rceil_{2} = 140, \left\lceil\!\!\left\lceil\!
		\begin{matrix}
			4 \\ 3
		\end{matrix}
		\!\right\rceil\!\!\right\rceil_{2}= 28, \left\lceil\!\!\left\lceil\!
		\begin{matrix}
			4 \\ 4
		\end{matrix}
		\!\right\rceil\!\!\right\rceil_{2} = 1$. 
		Therefore, we get
		\begin{align*}
			0  \cdot 1 + 100x + 140 (x^{2}-x) + 28 (x^{3} -5x^{2} +4x) + 1 (x^{4} - 14x^{3} + 49x^{2} -36x) 
		\end{align*}
		\begin{align*}
			&= 100x + 140x^{2} -140x + 28x^{3} -140x^{2} +112x + x^{4} -14x^{3} + 49x^{2} -36x \\ &= x^{4} + 14x^{3} + 49x^{2} + 36x \\ &= x(x+1)(x+4)(x+9) =  x^{\overline{\overline{4}}}_{2}.
		\end{align*}
	\end{ex}

	Let $\overline{\begin{NiceArrayWithDelims}{\llbracket}{\rrbracket}{c}
			n \\ k 
	\end{NiceArrayWithDelims}}_{s}$, $\overline{\left\{\mkern-4mu\braceVectorstack{n\\k}\mkern-4mu\right\}}_{s}$ and $\overline{\left\lceil\!\!\left\lceil\!
		\begin{matrix}
			n \\ k
		\end{matrix}
		\!\right\rceil\!\!\right\rceil}_{s}$ be the signed Stirling numbers with higher level and Lah numbers with order $s$, defined as
	
	\begin{align}
		\label{signedLah}
		\overline{\begin{NiceArrayWithDelims}{\llbracket}{\rrbracket}{c}
				n \\ k 
		\end{NiceArrayWithDelims}}_{s} = (-1)^{n+k} \begin{NiceArrayWithDelims}{\llbracket}{\rrbracket}{c}
			n \\ k 
		\end{NiceArrayWithDelims}_{s}  \nonumber \\
		\overline{\left\{\mkern-4mu\braceVectorstack{n\\k}\mkern-4mu\right\}}_{s} = (-1)^{n+k} \left\{\mkern-4mu\braceVectorstack{n\\k}\mkern-4mu\right\}_{s}  \nonumber \\
		\overline{\left\lceil\!\!\left\lceil\!
			\begin{matrix}
				n \\ k
			\end{matrix}
			\!\right\rceil\!\!\right\rceil}_{s} = (-1)^{n+k} \left\lceil\!\!\left\lceil\!
		\begin{matrix}
			n \\ k
		\end{matrix}
		\!\right\rceil\!\!\right\rceil_{s}.
	\end{align} The notation $\overline{\lah{n}{k}}$ was used for signed Lah numbers in \cite[pp. 66]{Mezo}.
	The equivalent inverse relation of (\ref{rising}) is given by the following theorem.
	
	\begin{thm}
		For $n \geq 0$ we have
		\begin{equation}
			\label{risinginverse}
			x_{s}^{\underline{\underline{n}}}= \sum_{k=0}^{n} \overline{\left\lceil\!\!\left\lceil\!
				\begin{matrix}
					n \\ k
				\end{matrix}
				\!\right\rceil\!\!\right\rceil}_{s} x_{s}^{\overline{\overline{k}}}.
		\end{equation}
	\end{thm}
	
	\begin{proof}
		By substituting $-x$ for $x$ in (\ref{rising}), and replacing $(-x)_{s}^{\underline{\underline{k}}}$ with $(-1)^{k} x_{s}^{\overline{\overline{k}}}$, since \[\prod_{i=0}^{k-1}(-x-i^{s}) = (-1)^{k} \prod_{i=0}^{k-1}(x+i^{s}),\] we get 
		\begin{align}
			(-x)^{\overline{\overline{n}}} &= \sum_{k=0}^{n} \left\lceil\!\!\left\lceil\!
			\begin{matrix}
				n \\ k
			\end{matrix}
			\!\right\rceil\!\!\right\rceil_{s} (-x)_{s}^{\underline{\underline{k}}} \nonumber \\
			&= \sum_{k=0}^{n} (-1)^{k} \left\lceil\!\!\left\lceil\!
			\begin{matrix}
				n \\ k
			\end{matrix}
			\!\right\rceil\!\!\right\rceil_{s}   x_{s}^{\overline{\overline{k}}}. \label{proofinverse}
		\end{align}
		Note that $(-1)^{n} (-x)_{s}^{\overline{\overline{n}}} = x_{s}^{\underline{\underline{n}}}$. Therefore, multiplying (\ref{proofinverse}) by $(-1)^{n}$, we get 
		\begin{align*}
			(-1)^{n} (-x)^{\overline{\overline{n}}} = \sum_{k=0}^{n} (-1)^{n+k} \left\lceil\!\!\left\lceil\!
			\begin{matrix}
				n \\ k
			\end{matrix}
			\!\right\rceil\!\!\right\rceil_{s}   x_{s}^{\overline{\overline{k}}}.
		\end{align*}
		The result follows.
	\end{proof}

	\begin{figure}[htbp] 
		\Large
		$$
		\xymatrix@C+2.5em@R+2.5em{
			&
			x^n
			\ar[ddr]_*{\large _{\large \left\{\mkern-8mu\braceVectorstack{n\\k}\mkern-8mu\right\}_{\Large s}}} \large
			\ar@/_15pt/[ddl]_*{\small \overline{\left\{\mkern-8mu\braceVectorstack{n\\k}\mkern-8mu\right\}}_{\large s}}   \huge
			& \\
			&  & \\
			x_{s}^{\overline{\overline{n}}}
			\ar[rr]^*{\large \left\lceil\!\!\left\lceil\!
				\begin{matrix}
					n \\ k
				\end{matrix}
				\!\right\rceil\!\!\right\rceil_{\large s}} \huge
			\ar[uur]_*{\large \begin{NiceArrayWithDelims}{\llbracket}{\rrbracket}{c}
					n \\ k 
				\end{NiceArrayWithDelims}_{\large s}} \huge
			& &
			x_{s}^{\underline{\underline{n}}} 
			\ar@/^15pt/[ll]^*{\large \overline{\left\lceil\!\!\left\lceil\!
					\begin{matrix}
						n \\ k
					\end{matrix}
					\!\right\rceil\!\!\right\rceil}_{\large s}} \huge
			\ar@/_15pt/[uul]_*{\large \overline{\begin{NiceArrayWithDelims}{\llbracket}{\rrbracket}{c}
						n \\ k 
				\end{NiceArrayWithDelims}}_{\large s}}
		}
		$$
		\caption{A diagram showing how Stirling numbers with higher level of both kinds and Lah numbers of order $s$ give coefficients for changing one basis of polynomials to another (the standard basis and the bases of falling and rising factorials with higher level).}
		\label{figure1}
	\end{figure}
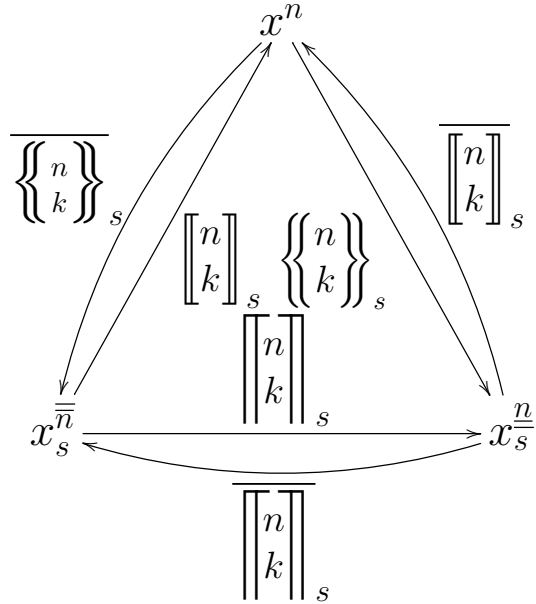

	\normalsize
	
	\subsection{Matrix representation}
	
	The Stirling numbers with higher level and Lah numbers of order $s$ satisfy the identities
	
	\begin{align}
		{x_{s}^{\overline{\overline{n}}}} = \sum_{k=0}^{n} \begin{NiceArrayWithDelims}{\llbracket}{\rrbracket}{c}
			n \\ k 
		\end{NiceArrayWithDelims}_{s} x ^{k} \label{matrix1} \\
		x^{n} = \sum_{k=0}^{n} \left\{\mkern-4mu\braceVectorstack{n\\k}\mkern-4mu\right\}_{s} {x_{s}^{\underline{\underline{n}}}}\label{matrix2} \\
		x^{\overline{\overline{n}}}_{s}= \sum_{k=0}^{n} \left\lceil\!\!\left\lceil\!
		\begin{matrix}
			n \\ k
		\end{matrix}
		\!\right\rceil\!\!\right\rceil_{s} x^{\underline{\underline{k}}}_{s}, \label{matrix3}
	\end{align} and the signed Stirling numbers with higher level and Lah numbers of order $s$ satisfy the identities
	
	\begin{align}
		{x_{s}^{\underline{\underline{n}}}} = \sum_{k=0}^{n} \overline{\begin{NiceArrayWithDelims}{\llbracket}{\rrbracket}{c}
				n \\ k 
		\end{NiceArrayWithDelims}}_{s} x ^{k} \label{matrix4} \\
		x^{n} = \sum_{k=0}^{n} \overline{\left\{\mkern-4mu\braceVectorstack{n\\k}\mkern-4mu\right\}}_{s} {x_{s}^{\overline{\overline{n}}}} \label{matrix5} \\
		x^{\underline{\underline{n}}}_{s}= \sum_{k=0}^{n} \overline{\left\lceil\!\!\left\lceil\! \label{matrix6}
			\begin{matrix}
				n \\ k
			\end{matrix}
			\!\right\rceil\!\!\right\rceil}_{s} x^{\overline{\overline{k}}}_{s}.
	\end{align}
	
	Let $\mathbb{Q} [x]$ denote the ring of polynomials with rational coefficients in the variable $x$. The ring  $\mathbb{Q} [x]$ can be seen as an infinite-dimensional vector space over the field $\mathbb{Q}$. In addition to the standard basis 
	
	\begin{align*}
		\Sigma = \{1, x, x^{2}, x^{3}, ...\}
	\end{align*} of the space $\mathbb{Q} [x]$, we have the natural bases 
	
	\begin{align*}
		\Sigma_{f} = \{1, x, x_{s}^{\underline{\underline{2}}}, x_{s}^{\underline{\underline{3}}}, ... \}, \\
		\Sigma_{r} = \{1, x, x_{s}^{\overline{\overline{2}}}, x_{s}^{\overline{\overline{3}}}, ... \}.
	\end{align*} Since all three of the above sets are bases of the space $\mathbb{Q} [x]$, every polynomial in each of them can be written as a linear combination of polynomials in any other of them. Identities (\ref{matrix1}) - (\ref{matrix3}) and (\ref{matrix4}) - (\ref{matrix6}) are infinite matrices \\

	$\left[ \begin{NiceArrayWithDelims}{\llbracket}{\rrbracket}{c}
		n \\ k 
	\end{NiceArrayWithDelims}_{s}\right]_{n, k=0}^{\infty}$, $\left[ \left\{\mkern-4mu\braceVectorstack{n\\k}\mkern-4mu\right\}_{s} \right]_{n, k=0}^{\infty}$,  $\left[\left\lceil\!\!\left\lceil\!
	\begin{matrix}
		n \\ k
	\end{matrix}
	\!\right\rceil\!\!\right\rceil_{s} \right]_{n, k=0}^{\infty}$ and \\
	
	$\left[\overline{\begin{NiceArrayWithDelims}{\llbracket}{\rrbracket}{c}
			n \\ k 
	\end{NiceArrayWithDelims}}_{s} \right]_{n, k=0}^{\infty}$,  $\left[ \overline{\left\{\mkern-4mu\braceVectorstack{n\\k}\mkern-4mu\right\}}_{s} \right]_{n, k=0}^{\infty}$,  $\left[ \overline{\left\lceil\!\!\left\lceil\!
		\begin{matrix}
			n \\ k
		\end{matrix}
		\!\right\rceil\!\!\right\rceil}_{s} \right]_{n, k=0}^{\infty}$, which are transition matrices between these bases.

	\begin{ex}
		Let us construct the standard and the natural basis for $\left[ \left\lceil\!\!\left\lceil\!
		\begin{matrix}
			n \\ k
		\end{matrix}
		\!\right\rceil\!\!\right\rceil_{2} \right]_{n, k=0}^{n, k = 7}$. We get
		\begin{align*}
			\begin{bmatrix}
				1 \\
				x  \\
				x_{2} ^ {\overline{\overline{2}}} \\
				x_{2} ^ {\overline{\overline{3}}} \\
				x_{2} ^ {\overline{\overline{4}}} \\
				x_{2}^ {\overline{\overline{5}}}\\
				x_{2} ^ {\overline{\overline{6}}} \\
			\end{bmatrix}
			= 
			\begin{bmatrix}
				1 & 0 & 0 & 0 &  0 & 0 & 0  \\
				0 & 1 & 0 & 0 & 0 & 0 & 0  \\
				0 & 2 & 1 & 0 & 0 &  0 & 0   \\
				0 & 10 & 10 & 1 & 0 & 0 & 0  \\
				0 & 100 & 140 & 28 & 1 & 0 & 0  \\
				0 & 1700 & 2900 & 840 & 60 & 1 & 0  \\
				0 & 44200 & 85800 & 31460 & 3300 & 110 & 1  \\
			\end{bmatrix}
			\begin{bmatrix}
				1 \\
				x\\
				x_{2}^ {\underline{\underline{2}}} \\
				x_{2} ^ {\underline{\underline{3}}} \\
				x_{2} ^ {\underline{\underline{4}}}\\
				x_{2} ^ {\underline{\underline{5}}} \\
				x_{2} ^{\underline{\underline{6}}} \\
			\end{bmatrix}.
		\end{align*}
		For example, $	x_{2} ^ {\overline{\overline{4}}} = 0 \cdot 1 + 100x + 140 	x_{2} ^ {\underline{\underline{2}}} + 28 	x_{2} ^ {\underline{\underline{3}}} + 1 \cdot 	x_{2} ^ {\underline{\underline{4}}}.$
	\end{ex}

	\subsection{The Lah Polynomials of Order $s$}
	
	We start by introducing the Lah polynomials of order $s$.
	
	\begin{defn}
		For $n \geq 1$, $s \geq 1$, the Lah polynomials of order $s$ are defined by
		\begin{align}
			\label{slahpoldef}
			\mathcal{L}_{n}^{s}(x) := \sum_{k=0}^{n} \left\lceil\!\!\left\lceil\!
			\begin{matrix}
				n \\ k
			\end{matrix}
			\!\right\rceil\!\!\right\rceil_{s} x^k.
		\end{align}
	\end{defn}
	
	\begin{thm}
		For $n \geq 1$, the Lah polynomials of order $s$ satisfy the recurrence relation
		\begin{align}
			\label{slahrecpol}
			\mathcal{L}_{n+1}^{s}(x) = x \mathcal{L}_{n}^{s}(x) + n^{s} \mathcal{L}_{n}^{s}(x) + B_{n}^{s}(x),
		\end{align}
		where $B_{n}^{s}(x)$ is for $n \geq 1$, $s \geq 1$ defined by
		\begin{align*}
			B_{n}^{s}(x) = \sum_{k=0}^{n} k^{s} \left\lceil\!\!\left\lceil\!
			\begin{matrix}
				n \\ k
			\end{matrix}
			\!\right\rceil\!\!\right\rceil_{s} x^{k}. 
		\end{align*}
	\end{thm}
	
	\begin{proof}
		Using recurrence relation (\ref{slahrec}) and definition (\ref{slahpoldef}), we get 
		\begin{align*}
			\mathcal{L}_{n+1}^{s}(x) &= \sum_{k=0}^{n+1} \left\lceil\!\!\left\lceil\!
			\begin{matrix}
				n+1 \\ k
			\end{matrix}
			\!\right\rceil\!\!\right\rceil_{s} x^{k} \\
			&= \left\lceil\!\!\left\lceil\!
			\begin{matrix}
				n+1 \\ 0
			\end{matrix}
			\!\right\rceil\!\!\right\rceil_{s} + \sum_{k=1}^{n} \left\lceil\!\!\left\lceil\!
			\begin{matrix}
				n+1 \\ k
			\end{matrix}
			\!\right\rceil\!\!\right\rceil_{s} x^{k} + \left\lceil\!\!\left\lceil\!
			\begin{matrix}
				n+1 \\ n+1
			\end{matrix}
			\!\right\rceil\!\!\right\rceil_{s} x^{n+1} \\
			&= \delta_{n+1,0}  + \sum_{k=1}^{n} \left( \left\lceil\!\!\left\lceil\!
			\begin{matrix}
				n \\ k-1
			\end{matrix}
			\!\right\rceil\!\!\right\rceil_{s} + (n^{s} +k^{s}) \left\lceil\!\!\left\lceil\!
			\begin{matrix}
				n \\ k
			\end{matrix}
			\!\right\rceil\!\!\right\rceil_{s} \right) x^{k} + x^{n+1} \\
			&= \sum_{k=0}^{n-1} \left\lceil\!\!\left\lceil\!
			\begin{matrix}
				n \\ k
			\end{matrix}
			\!\right\rceil\!\!\right\rceil_{s} x^{k+1} + x^{n+1} + n^{s} \sum_{k=1}^{n} \left\lceil\!\!\left\lceil\!
			\begin{matrix}
				n \\ k
			\end{matrix}
			\!\right\rceil\!\!\right\rceil_{s} x^{k} + \delta_{n+1,0} + \sum_{k=1}^{n} k^{s} \left\lceil\!\!\left\lceil\!
			\begin{matrix}
				n \\ k
			\end{matrix}
			\!\right\rceil\!\!\right\rceil_{s} x^{k} \\
			&= x \sum_{k=0}^{n} \left\lceil\!\!\left\lceil\!
			\begin{matrix}
				n \\ k
			\end{matrix}
			\!\right\rceil\!\!\right\rceil_{s} x^{k} + n^{s} \sum_{k=0}^{n} \left\lceil\!\!\left\lceil\!
			\begin{matrix}
				n \\ k
			\end{matrix}
			\!\right\rceil\!\!\right\rceil_{s} x^{k} + \sum_{k=1}^{n} k^{s} \left\lceil\!\!\left\lceil\!
			\begin{matrix}
				n \\ k
			\end{matrix}
			\!\right\rceil\!\!\right\rceil_{s} x^{k} \\
			&= x\mathcal{L}_{n}^{s}(x) + n^{s} \mathcal{L}_{n}^{s}(x) + B_{n}^{s}(x).
		\end{align*}
		Note that $\delta_{n+1,0}=0$ in our case.
	\end{proof}
	
	\begin{rem}
		The definition of the polynomials $B_{n}^{s}(x)$ is equivalent to 
		\begin{align*}
			B_{n}^{s}(x)=  \sum_{k=1}^{n} k^{s} \left\lceil\!\!\left\lceil\!
			\begin{matrix}
				n \\ k
			\end{matrix}
			\!\right\rceil\!\!\right\rceil_{s} x^{k},
		\end{align*} since the term corresponding to $k=0$ equals $0$. 
	\end{rem}
	
	\begin{lem}
		\label{lemma}
		For all $n \geq 1$, $s \geq 1$, the polynomials $B_{n}^{s} (x)$ satisfy
		\begin{align*}
			B_{n}^{s} (x) = \sum_{i=0}^{s} \sttwo{s}{i}  x^{i} \frac{d^{i}}{dx^{i}} \mathcal{L}_{n}^{s} (x).
		\end{align*}
	\end{lem}
	\begin{proof}
		It follows from (\ref{stirpolid}) that 
		\begin{align}
			B_{n}^{s}(x) &= \sum_{k=0}^{n} \left\lceil\!\!\left\lceil\!
			\begin{matrix}
				n \\ k
			\end{matrix}
			\!\right\rceil\!\!\right\rceil_{s} k^{s} x^{k} \nonumber \\
			&= \sum_{k=0}^{n} \left\lceil\!\!\left\lceil\!
			\begin{matrix}
				n \\ k
			\end{matrix}
			\!\right\rceil\!\!\right\rceil_{s} \sum_{i=0}^{s} \sttwo{s}{i} k^{\underline{i}} x^{k} \nonumber \\
			&= \sum_{i=0}^{s} \sttwo{s}{i} x^{i} \sum_{k=0}^{n} k^{\underline{i}} \left\lceil\!\!\left\lceil\!
			\begin{matrix}
				n \\ k
			\end{matrix}
			\!\right\rceil\!\!\right\rceil_{s}  x^{k-i} \nonumber \\
			&= \sum_{i=0}^{s} \sttwo{s}{i} x^{i} \frac{d^{i}}{dx^{i}} \mathcal{L}_{n}^{s} (x). \label{bnsx}
		\end{align}
	\end{proof} Since $s \geq 1$ by the definition of $\mathcal{L}_{n}^{s}(x)$ and $B_{n}^{s}(x)$ (the term corresponding to $i=0$ in (\ref{bnsx}) is nonzero only when $i=s=0$), (\ref{slahrecpol}) can be written as
	
	\begin{align*}
		\mathcal{L}_{n+1}^{s}(x) = x \mathcal{L}_{n}^{s}(x) + n^{s} \mathcal{L}_{n}^{s}(x) + \sum_{i=1}^{s} \sttwo{s}{i} x^{i} \frac{d^{i}}{dx^{i}} \mathcal{L}_{n}^{s} (x).
	\end{align*}

	\section{Connections Between the Lah Numbers with Higher Level and the Lah Numbers of Order $s$}
	
	Here, we give some connections between the Lah numbers with higher level and the Lah numbers of order $s$. Indeed, an obvious one is 
	
	\begin{align*}
		\left\lfloor\!\!\left\lfloor\!
		\begin{matrix}
			n \\ k
		\end{matrix}
		\!\right\rfloor\!\!\right\rfloor_{1} =
		\left\lceil\!\!\left\lceil\!
		\begin{matrix}
			n \\ k
		\end{matrix}
		\!\right\rceil\!\!\right\rceil_{1} = \lah{n}{k}.
	\end{align*}
	
	Next we notice the following.
	
	\begin{prop}
		\label{inequality1}
		For all $n, k, s \in \mathbb{N}$, we have the inequality
		\begin{align*}	
			\left\lfloor\!\!\left\lfloor\!
			\begin{matrix}
				n \\ k
			\end{matrix}
			\!\right\rfloor\!\!\right\rfloor_{s} \geq 
			\left\lceil\!\!\left\lceil\!
			\begin{matrix}
				n \\ k
			\end{matrix}
			\!\right\rceil\!\!\right\rceil_{s}.
		\end{align*}
	\end{prop}
	
	\begin{proof}
		Applying the recurrence relations (\ref{slah}) and (\ref{slahrec}), we need to show that $(n+k-1)^{s} \geq (n-1)^{s} + k^{s}$ for all $n, k, s \in \mathbb{N}$. Let $a=n-1$ and $b=k$. We need to prove that $(a+b)^{s} \geq a^{s} + b^{s}$. From the binomial theorem we observe that two of the terms in the expansion of $(a+b)^{s}$ are $a^{s}$ and $b^{s}$, while all other terms are non-negative. 
	\end{proof}
	
	An extended version of the previous theorem is the following inequality.
	
	\begin{thm}
		\label{inequality2}
		For all $n, k, s \in \mathbb{N}$, we have the inequality
		\begin{align*}
			\left\lfloor\!\!\left\lfloor\!
			\begin{matrix}
				n \\ k
			\end{matrix}
			\!\right\rfloor\!\!\right\rfloor_{s} \geq 
			\left\lceil\!\!\left\lceil\!
			\begin{matrix}
				n \\ k
			\end{matrix}
			\!\right\rceil\!\!\right\rceil_{s} \geq \begin{NiceArrayWithDelims}{\llbracket}{\rrbracket}{c}
				n \\ k
			\end{NiceArrayWithDelims}_{s} \geq \left\{\mkern-4mu\braceVectorstack{n\\ k} \mkern-4mu\right\}_{s}.
		\end{align*}
	\end{thm}
	
	\begin{proof}
		Applying recurrence relations  (\ref{sstirling1rec}), (\ref{sstirling2rec}), (\ref{slah}) and (\ref{slahrec}), we need to show that $(n+k-1)^{s} \geq (n-1)^{s} + k^{s} \geq (n-1)^{s} \geq k^{s}$. In the proof of Proposition \ref{inequality1}, we already proved that $(n+k-1)^{s} \geq (n-1)^{s} + k^{s}.$ The inequality $(n-1)^{s} + k^{s} \geq (n-1)^{s}$ is obvious. The inequality $(n-1)^{s} \geq k^{s}$ for Stirling numbers of both kinds with higher level obviously holds for $n-1 \geq k$ and if $k \geq n$, $\begin{NiceArrayWithDelims}{\llbracket}{\rrbracket}{c}
			n \\ k
		\end{NiceArrayWithDelims}_{s} = \left\{\mkern-4mu\braceVectorstack{n\\ k} \mkern-4mu\right\}_{s}$ according to (\ref{cases1a}) and (\ref{n<k}). Therefore, $\begin{NiceArrayWithDelims}{\llbracket}{\rrbracket}{c}
			n \\ k
		\end{NiceArrayWithDelims}_{s} \geq \left\{\mkern-4mu\braceVectorstack{n\\ k} \mkern-4mu\right\}_{s}$. 
	\end{proof}

	\begin{thm}
		For $k=n-1$, we have
		\begin{align*}
			\left\lceil\!\!\left\lceil\!
			\begin{matrix}
				n \\ n-1
			\end{matrix}
			\!\right\rceil\!\!\right\rceil_{s} = \frac{
				\left\lfloor\!\!\left\lfloor\!
				\begin{matrix}
					n \\ n-1
				\end{matrix}
				\!\right\rfloor\!\!\right\rfloor_{s}}{2^{s-1}}.
		\end{align*}
	\end{thm}
	
	\begin{proof}
		Considering (\ref{k=n-1(1)}) and (\ref{k=n-1(2)}), we get
		\begin{align*}
			\label{connection1proof}
			\frac{\sum_{j=1}^{n-1}(2j)^{s}}{2^{s-1}} &\left.= \right. \sum_{j=1}^{n-1}2j^{s}.
		\end{align*}
		The result follows.
	\end{proof}

	\section*{Acknowledgement}
	The author is grateful to his parents for their support. The author is grateful to Professor Stephan Wagner for review of the article and very helpful comments and suggestions. Professor Wagner suggested the author to study the polynomials $Q_{n}^{s}(x)$, studied in Section \ref{qnsxsection}.

\end{document}